\documentclass[12pt]{amsart}
\usepackage{a4,verbatim}
\usepackage{amssymb}
\usepackage{amsmath}
\usepackage{amsthm}
\usepackage{amstext}
\usepackage{amscd}
\usepackage{latexsym}
\usepackage{graphics}
\usepackage{color}

\swapnumbers
\theoremstyle{plain}
\newtheorem{thm}{Theorem}[section]

\newtheorem{lemma}[thm]{Lemma}

\newtheorem{prop}[thm]{Proposition}
\newtheorem{conj}[thm]{Conjecture}
\theoremstyle{definition}

\newtheorem{rmk}[thm]{Remark}

\newcommand{\supp}{{\rm Supp}}

\renewcommand{\tilde}{\widetilde}

\newcommand{\C}{{\mathbb C}}

\renewcommand{\P}{{\mathbb P}}
\newcommand{\Q}{{\mathbb Q}}
\newcommand{\R}{{\mathbb R}}

\newcommand{\V}{{\mathbb V}}

\newcommand{\Z}{{\mathbb Z}}

\newcommand{\la}{\longrightarrow}
\newcommand{\round}[1]{\lfloor #1 \rfloor}

\input{xy}
\xyoption{all}

\begin{document}

\title{On Effective log Iitaka Fibration for 3-folds and 4-folds}
\author{Gueorgui Todorov and Chenyang Xu}

\begin{abstract}We prove the effectiveness of the log Iitaka fibration in Kodaira codimension two for varieties of dimension$\le 4$. In particular, we finish the proof of effective log Iitaka fibration in dimension two. Also, we show that for the log Iitaka fibration, if the fiber is of dimension two, the denominator of the moduli part is bounded.\end{abstract}

\maketitle


\vspace*{6pt}
\section{Introduction}
One of the main problems in birational 
algebraic geometry is to understand the 
structure of pluricanonical  maps.
If $X$ is a smooth complex projective variety  of dimension $n$, then 
we define the
pluricanonical maps $$\phi_{rK_X}:X\dashrightarrow \mathbb{P}(H^0(X,\mathcal{O}_X(rK_X)))$$
determined  by the linear system $|rK_X|$. If this linear system is nonempty for some natural number $r$, which is conjectually equivalent to $X$ being non-uniruled, then for $r$ sufficiently divisible the maps $\phi_{rK_X}$ become birational to a fixed algebraic fiber space $\phi:X'\la Y'$  called the Iitaka fibration of $X$ and the Kodaira dimension $\kappa(X)$ is set to be $\dim Y'$.
 It is then natural to look for a uniform $r$ (that is an integer $r$ which depends only on the dimension of $X$) for which we are inducing  a map birational to the Iitaka fibration. After the monumental work \cite{bchm} of proving the finite generation of the pluricanonical rings, the effective Iitaka fibration problem is largely related to finding generators which are in less than a uniform degree for a given dimension. If $\kappa(X)=\dim X$, then $X$ is called of general type and in this case
 Hacon-M$^\textrm{c}$Kernan \cite{hm06} and  Takayama \cite{tak06} following ideas of Tsuji, have shown that  for a smooth 
projective variety of general type and dimension $n$ there exists an integer 
$r_n$, which depends only on $n$, such that $\phi_{rK_X}$ is birational 
for $r \ge r_n$. The question remains widely open when the Kodaira dimension in not maximal and it is completely known only  up to dimension three 
(cf. \cite{fm,vz,rin}).

The current approach for studying this problem is to use Kawamata's canonical bundle formula which identifies the pluricanonical ring of $X$ with the pluricanonical ring of a pair $(Y,\Delta+L)$ of log general type (cf. \cite{fm}, \cite{ps09}),  where $(Y,\Delta)$ is a $klt$ pair  and $L$ is a $\Q$-line bundle coming from variation of Hodge structure. This raises the natural question of a log analogue of the above statement. Another setting where boundary divisors naturally occur is in moduli problems.
Because of the presence of a boundary divisor the best statement that we can hope for is the following: for $(X,D)$ a $klt$ pair of Kodaira dimension $\kappa \ge0$ there is a natural number $m$, that depends only on the coefficients of $\Delta$  and the dimension of $X$ such that $|\round{m(K_X+\Delta)}|$ is birational to the Iitaka fibration.
To be able to  use the canonical bundle formula inductively we have to allow that the coefficients of $\Delta$ lie in a possibly infinite set of rational numbers.  It turns out that the natural assumption on the set of coefficients is the {\it descending chain condition} (DCC). At the first glance, it may seem to be quite technical and somewhat artificial. However, this is  a very natural and useful condition arising in many problems. The interested reader is recommended to read \cite{sh92}, \cite{ko94} for more detailed discussions of this. 
Thus the general conjecture is formulated as follows:
\begin{conj}(Effective Log Iitaka Fibration Conjecture)
 Let $(X,\Delta)$ be a $klt$ (or $\epsilon-lc$) pair such that $K_X+\Delta$ is pseudo-effective and the coefficients in $\Delta$ are in a DCC set $\mathcal{A}$. Then there is a constant  $r$ depending only on the dimension $X$ and on $\mathcal{A}$ such that $|\round{r(K_X+\Delta)}|$ induces a birational map. 
\end{conj}

In this paper we prove the following:
\begin{thm}\label{THM} Let $(X,\Delta)$ be a klt pair of dimension two,
 three or four and Kodaira codimension two. Assume that the coefficients of $\Delta$ are in a 
DCC set of rational numbers $\mathcal{A}\subset [0,1]$.  Then there is an explicitly computable
 constant $m$ depending only on 
the set $\mathcal{A}$  such that
$\round{ m(K_X+\Delta) }$ induces the Iitaka fibration.
\end{thm} 

In particular, this completes the effectiveness of log Iitaka fibration in dimension two (for the other cases, see \cite{tod}).
The dimension two case  relies on results by Alexeev and Alexeev-Mori \cite{alex,am04}. 

A key tool in the dimension three and the dimension four case is the canonical bundle formula (cf. \cite{ka98}, \cite{fm}, \cite{ko07}). In our case 
it roughly says that if $f:X\la Y$ is the Iitaka fibration for $K_X+\Delta$ and  $K_X+\Delta=f^*D$ for some $\Q$-divisor $D$
 on $Y$ then we can define the  \emph{discriminant}  or {\it divisorial part} on $Y$ for  $K_X+\Delta$ to be the  $\Q$-Weil divisor $B_Y:=\sum_P b_P P$, where $1-b_P$ is the maximal
 real number $t$ such that the log pair $(X,\Delta+tf^*(P))$ has log canonical singularities over the
 generic point of $P$. The sum runs over all codimension one points of $Y$, but it has finite
support. The \emph{moduli part} or \emph{J-part} is the $\Q$-line bundle $M_Y$ on $Y$ which satisfies
$$ K_X+\Delta= f^*(K_Y+B_Y+M_Y). $$

If $F$ is the general fiber of $f$ and if $h^0(F,b(K_F+\Delta_{|F}))\neq 0$, then for every integer $r$ divisible by $b$ we have 
$$H^0(X,\round{r(K_X+\Delta)})=H^0(Y,\round{r(K_Y+B_Y+M_Y)}).$$

Using ideas of Mori and Fujino,  and the dimension two case of Theorem \ref{THM}, we prove an analogue of Mori and Fujino bounding the denominators in relative dimension two.
\begin{thm}\label{rel2} When the relative dimension of $f$ is two, there is a natural number $m$ depending only on the coefficients of $\Delta$ such that $mM$ is Cartier.
\end{thm}

 In the  dimension three case of Theorem \ref{THM}  the base is a curve and the conclusions  follows (\ref{rel2}) by standard arguments as in \cite{fm}. When the dimension is four a key
 observation is that many results in \cite{am04} can be improved by 
adding a nef divisor (\ref{neftail}). 
 Similar statements were  essentially proved by Viehweg-Zhang \cite{vz}, but we will simplify their proof by putting it in the context of \cite{alex}, \cite{am04} and the recent main stream investigations on adjoint linear systems. We should also  remark that in (\cite{ps09}, Conjecture 7.13), Prokhorov-Shokurov list a series of conjectures concerning  the solution of the effective Iitaka fibration problem. Putting our work in the frame of their conjectures, we indeed show Conjecture 7.13(2) in the relative dimension two case. However, we do not prove the semi-ampleness statements of Conjecture 7.13 (1) and (3).

The paper is structured as follow. In Section 2, we show that the DCC assumption on the coefficients of the boundary indeed forces  that the $klt$ surface pairs of log Kodaira dimension zero to be $\epsilon$-log canonical, for some $\epsilon$ depending on the set of coefficients. In Section 3, we prove that under the same assumption there is a uniform $b$ such that $b(K_S+B)\sim 0$. In Section 4, using ideas of Mori and Fujino, we deduce Theorem (\ref{rel2}). And then we complete the proof of Theorem (\ref{THM}).

{\bf Acknowledgment:} The authors would like to thank Christopher Hacon,  J\'anos Koll\'ar and James M$^\textrm{c}$Kernan for numerous conversations and suggestions. This material is also based upon work when the second named author is in Institute for Advanced Study and supported by the NSF under agreement No. DMS-0635607. Any opinions, findings and conclusions or recommendations expressed in this
material are those of the author and do not necessarily reflect the views of the
NSF.

\subsection{Notations and Conventions.}

We will work over the field of complex numbers $\mathbb{C}$. 
A $\Q$-Cartier divisor $D$ is nef if $D\cdot C \ge0$ for any curve $C$ on $X$.
 We call two $\Q$-divisors $D_1, D_2$ $\Q$-linearly equivalent $D_1\sim_\Q D_2$ if 
there exists an integer $m>0$ such that $mD_i$ are integral and linearly equivalent.
We call two $\Q$-Cartier divisors $D_1, D_2$ numerically equivalent $D_1\equiv D_2$ if $(D_1-D_2)\cdot C=0$  for any curve $C$ on $X$.
A log pair $(X,\Delta)$ is a normal variety $X$ and an effective $\Q$-Weil 
divisor $\Delta$ such that $K_X+\Delta$ is $\Q$-Cartier. A projective morphism 
$\mu:Y \la X$ is a log resolution of the pair  $(X,\Delta)$ if $Y$ is smooth and 
$\mu^{-1}(\Delta)\cup\{\textrm{exceptional set of } \mu\}$ is a divisor 
with simple normal crossing support. For such $\mu$ we write   
$\mu^*(K_X+\Delta)  =K_Y+\Gamma$, and $\Gamma=\Sigma a_i\Gamma_i$ 
where $\Gamma_i$ are distinct integral divisors. A pair is called 
$klt$ (resp. $lc$, $\epsilon$-$lc$) if there is a log resolution $\mu:Y \la X$ such that in 
the above notation we have $a_i <1$ (resp. $a_i\le 1,a_i \le 1-\epsilon$). The number $1-a_i$ is called 
log discrepancy of $\Gamma_i$ with respect to the pair  $(X,\Delta)$.

For conventions and results about DCC set we refer to  \cite{am04}, in particular, (\cite{am04}, 2.2-2.7 and 3.4-3.6). When we say some quantity is bounded, it always means there is a computable bound depending on the data we give. We will not keep tracking of the explicit bound, however, this does not require much effort by following the arguments in \cite{am04}.  Most of the results about surfaces  that we use appear already in \cite{alex}, but to make the constants explicitly computable we refer to the more recent  \cite{am04}. 

\section{$\epsilon$-log canonicity}

The main result of this section is to show the following
\begin{thm}\label{epsilon}
Let $(S,B=\sum b_iB_i)$ be a $klt$ projective surface pair. If $K_S+B\equiv 0$ and the coefficients of $B$ are in a DCC set $\mathcal{A}$, then there exists an $\epsilon=\epsilon(\mathcal{A})>0$, which only depends on $\mathcal{A}$ such that  $(S,B)$ is $\epsilon$-log canonical. 
\end{thm}
\begin{proof} We can run a minimal model program for $S$,
 $$S=S_0\to S_1 \to \cdots \to S_n.$$
By pushing  forward $B$ to $S_n$, we see that  it suffices to prove the statement for $S_n$. So we need only to prove for the case that $S$ is a log del Pezzo surface of Picard number 1, $S$ admits a Fano contraction to some curve with a general fiber $\P^1$ or $S$ has $K_S \equiv 0$. This is done in (\ref{-2}), (\ref{-1}) and (\ref{0}).

\end{proof}

The main tool we use is the following theorem,
\begin{thm}\cite[Theorem 3.2]{am04}\label{dp}
Let $\mathcal{A}\subset [0,1]$ be an arbitrary DCC set. Then there exists a constant $\delta=\beta(\mathcal{A})$ depending only on $\mathcal{A}$ so that the following holds. Let $X$ be a normal projective surface, $B_j$ be divisors on $X$, and let $b_j$, $x_j$ be positive real numbers. Assume that
\begin{enumerate}
 \item $X$ is a singular $\Q$-factorial surface with $-K_X$ ample and $\rho(X)=1$,
\item $b_j>0$ and $b_j \in \mathcal{A}$,
\item $1-\delta<x_j\le 1$,
\item at least one $x_j$ is strictly less than 1,
\item the pair $(X,\sum x_jb_jB_j)$ is lc. 
\end{enumerate}
Then the divisor $K_X+\sum x_jb_jB_j$ is not numerical equivalent to 0.
 
\end{thm}

Also we need the following lemma on DCC sets,

\begin{lemma}\cite[Lemma 3.6]{am04}\label{dcc}
 Let $\mathcal{A}$ be a DCC set containing 1. Let $a,\delta>0$ such that $\delta\le mf_2(\mathcal{A},a)$ (see \cite{am04} (3.5) for the definition of the function $mf_2(\mathcal{A},a)$). Then for any finite number of arbitrary $b_j$, $x_j\in R$ such that $0<b_j\in \mathcal{A}$ and $1-\delta<x_j\le 1$ for all $j$, $ x_j<1$ for some $j$, then we have $\sum x_jb_j\neq a$.
\end{lemma}

\begin{lemma}\label{-2}
 If $S$ is a log del Pezzo surface of Picard number 1. $B=\sum b_jB_j$ is a $\Q$-divisor such that $K_S+B\equiv 0$, $(S,B)$ is $klt$ and the coefficients of $B$ are in a DCC set $\mathcal{A}$. Then $S$ is $\epsilon=\epsilon(\mathcal{A})$-log canonical. 
\end{lemma}

\begin{proof}
We assume $1\in \mathcal{A}$. If $(S,B)$ is not $\epsilon$-log canonical, then we can extract a divisor $f: S'\to S$ such that 
$$f^*(K_S+\sum b_jB_j)\equiv K_{S'}+\sum b_jB_j+(1-a)E,$$
with $a< \epsilon$. $S'$ has Picard number 2, so it has two extremal rays, one
of which  is given by $E$. Contracting the other extremal ray we obtain a morphism $S'\to S_0$ to a surface $S_0$. In fact if $S_0$ were a curve,  by  taking the intersection of $K_{S'}+\sum b_jB_j+(1-a)E$ with a general fiber, we have
$$\sum n_ib_j+n(1-a)=2,$$  
where $n_i\ge 0$  and  $n>0$ are integers, which contradicts  (\ref{dcc}), provided we choose $\epsilon=mf_2(\mathcal{A},2)$.
 Then $S_0$ is a log del Pezzo of Picard number 1.
Now  $g$ does not contract $E$, and 
$$g_*(K_{S'}+\sum b_jB_j+(1-a)E)=K_{S_0}+(1-a)E+\sum b_jg_*(B_j)\equiv 0.$$ 
Since $K_S+B\equiv 0$    we still have that 
$(S_0,(1-a)E+\sum b_jg_*(B_j))$ is a $klt$ pair. Then Theorem \ref{dp} applied 
 to $S_0$ with $b_E=1, x_E=(1-a)$, all other $x_j=1$ and  $\epsilon=\delta$ implies that $K_{S_0}+(1-a)E+\sum b_jg_*(B_j)$ is not numerically equivalent to zero which is a contradiction. 
 \end{proof}

\begin{lemma}\label{-1}
Let $(S,B)\to C$ be  a morphism from a surface to a curve, such that the generic fiber is $\P^1$ and let $(S,B)$ be $klt$ with the coefficients of $B$  in a DCC set $\mathcal{A}$. Then $S$ is $\epsilon=\epsilon(\mathcal{A})$-log canonical.   
\end{lemma}
\begin{proof}
 Suppose that this is not the case. Then extract a divisor $f: S'\to S$ such that 
$$f^*(K_S+\sum b_jB_j)\equiv K_{S'}+\sum b_jB_j+(1-a)E,$$
with $a< \epsilon$. Running a minimal model program for $(S',\sum b_jB_j)$ that does not contracts $E$, we end with a surface $S''$ that is either log del Pezzo or admits a Fano contraction to a curve. Furthermore the  coefficients of  $E$  in the pair $(S'',B'')$ is $(1-a)$ and $S''+B''\equiv 0$.

In the case $S''$ is log del Pezzo, we get a contradiction as in (\ref{-2}). 
If $S''$ admits a Fano contraction to a curve, by intersecting with a general fiber, we can apply the argument in (\ref{-2}) again.

\end{proof}

Finally we deal with the case in which $K_S$ is numerically trivial.

\begin{lemma}\label{0}
 There is $\epsilon>0$  such that every $klt$ surface with $K_S\equiv 0$ is $\epsilon$-log canonical. 
\end{lemma}
\begin{proof}
   Set $\epsilon=\min\{\beta(\{1\}), mf_2(\{1\}',2)\}$. Suppose that $S$ is not
$\epsilon$-log canonical and  extract a  divisor $E$ by $f:S'\to S$ such that
$$f^*(K_S)=K_{S'}+(1-a)E \mbox{ and } 0<a<\epsilon.$$ 
 Then running a minimal model for $S'$, which does not contract $E$,  we can apply one of the above two cases.

\end{proof}

\section{Bounding  the index}

Observe that on the varying  $a,b$ and $c$   the weighted projective spaces $\P(a,b,c)$ with three lines $(x_i=0)$ forms  an unbounded family of $lc$ surface pairs with coefficients are in the DCC set $\{1\}$. However by restricting to the $klt$ case we have the following boundedness result:

\begin{thm}\label{uniform}
 For $(S,B)$  a $klt$ pair  such that  $K_S+B\equiv 0$ and the coefficients of $B$ are in a DCC set $\mathcal{A}$,  there is a natural number $b=b(\mathcal{A})$, such that $b(K_S+B)$ is Cartier and $H^0(S,b(K_S+B))=1$.
\end{thm}

Form Theorem (\ref{epsilon}) we know that there is an $\epsilon>0$ such that $(S,B)$ is $\epsilon$-log canonical and we will assume  $0<\epsilon <1/\sqrt{3}$.

\begin{lemma}\cite[Lemma 1.2, Theorem 1.8]{am04}\label{picard}
Let $X$ be a non singular projective surface and $B=\sum b_jB_j$ be an $\R$-divisor on $X$ with $0\le b_j \le 1-\epsilon<1$. Assume $K_X+B\equiv 0$, then
\begin{enumerate}
\item  If $E$ is an irreducible curve on $X$ and $E^2<0$ then $E\cong \P^1$ and $E^2>-2/\epsilon$; and
\item  $\rho (X)\le 128/\epsilon^5$. 
 
\end{enumerate}

\end{lemma}

\begin{prop}\label{weil}
 $(S,B)$ as in (\ref{uniform}). There exists an integer $t=t(\epsilon)$ such that for any Weil divisor $D$ on $S$, $tD$ is Cartier.
\end{prop}

\begin{proof}
The argument is the parallel to  the one given in (\cite{am04}, Lemma 3.7), though there $S$ is assumed to be a del Pezzo surface of Picard number 1.

Take the minimal resolution $f:\tilde{S}\to S$ of $S$ and write 
$$f^*(K_S)=K_{\tilde{S}}+\sum a_iE_i.$$
Applying (\ref{picard}), we conclude that the determinant $t$ of the intersection matrix of the exceptional curves $(-F_iF_k)_{ik}$ is bounded by 
$$t\le [2/\epsilon]^{[128/\epsilon^5]}.$$
Since $S$ has  rational singularities, for any Weil divisor $D$ on $S$, $tD$ is a Cartier divisor.
\end{proof}

In the next Proposition we repeatedly use the trivial property of DCC set stated in the following Lemma.

\begin{lemma}\label{finiteDCC} Let $c$ be a positive real number and  $\mathcal{A}\subset [0,1]$ a DCC set. Then there are a finite number of ways to write $c$ as the sum of $a_i\in\mathcal{A}$.
\end{lemma}

\begin{prop}\label{index}
 
 $(S,B)$ as in (\ref{uniform}). There exists an integer $t=t(\epsilon)$ such that $t(K_S+B)$ is Cartier.

\end{prop}
\begin{proof}
Running a minimal model program for $S$, and we end with a $klt$ surface $f:S\to S'$. If $t(K_S'+B')$ is Cartier, since $$f^*(t(K_{S'}+B'))\sim t(K_S+f_*^{-1}B'+a_iE_i)$$ and
$$K_S+f_*^{-1}B'+a_iE_i\equiv K_S+B\equiv 0,$$
we know $K_S+f_*^{-1}B'+a_iE_i= K_S+B$, thus it suffices to prove for $(S',B')$. $S'$ is either of Kodaira dimension 0, has a contraction to a curve with general fiber $\P^1$ or a log del Pezzo surface, so we only need to prove the proposition in these cases.

{\it (Case 1)}:  $S$ has Kodaira dimension 0. Then $B=0$ and it follows directly from (\ref{weil});

{\it (Case 2)}: $S$ admits a contraction to a curve $f:S\to C$ with  general fiber $\P^1$. By taking the intersection of $(S,B)$ with the general fiber, we have 
$$2=\sum_{f(B_j)\neq pt} a_j, \mbox{ where } a_j \in \mathcal{A}'.$$
On the other hand, applying the canonical bundle formula in this simple case, we have 
$$K_S+B\equiv f^*(K_C+B_C)\equiv 0, \mbox{ where }B_C=\sum_{f(B_j)=P_j}\frac{a_j+n-1}{n}P_j $$
From Lemma (\ref{finiteDCC}) we know that there are only finitely many choices of $a_j$, and we choose a natural number $N=N(\mathcal{A})$ that clears  all the denominators of $a_j$.

{\it (Case 3)}:  $S$ is a log del Pezzo surface of Picard number one. For 
$t$ as in (\ref{weil}), $K_S+\sum b_jB_j\equiv 0$ implies that 
$$-tK_S^2=\sum tb_j B_j\cdot K_S.$$ 
Here $tK_S^2$ and $tB_j\cdot K_S$ are  integers. 
From the proof of (\cite{am04}, Lemma 3.7), we can bound $tK_S^2$ by
$$tK_S^2 \le [2/\epsilon]^{[128/\epsilon^5]} \times ([2/\epsilon]+2)^2.$$
Lemma \ref{finiteDCC} implies that  there are only finitely many possible $b_i$ appearing as the coefficients of the boundary divisor. So there is a positive number $N=N(\mathcal{A})$ such that $NK_S$ and $N(b_jB_j)$ are all Cartier divisors.

\end{proof}

Because of the above discussion, if we take the minimal resolution $\pi:\tilde{S}\to S$ of $S$, pull back $K_S+B$, then the denominators can be killed by a uniform multiple $t=t(\epsilon)$. Thus by running a minimal model program for $\tilde{S}$ we only need to prove (\ref{uniform}) in the case that $S$ is a smooth minimal surface and that the coefficients of $B$ all have the form $r/t$ for some uniform $t$. When $K_S\equiv 0$, and hence $B=0$ the result follows  from the classification theory of smooth surfaces; when $S$ is rational, a Cartier divisor is numerically trivial if and only if it is a trivial divisor. Thus the only tricky case is dealt with in the next Lemma.

\begin{lemma}
If $S$ is a smooth minimal ruled surface over a curve of positive genus and there is a $\Q$-divisor $B$ such that $(S,B)$ is $klt$ and $K_S+B\equiv 0$, then $S=E\times \P^1$, where $E$ is an elliptic curve. In particular, if for some integer $N$, $NB$ is an integer divisor, then $\mathcal{O}(N(K_S+B))\cong \mathcal{O}_S$.
\end{lemma}

\begin{proof} 
We use conventions for ruled surfaces as in  (\cite{ha77}, V.2). Let $\pi:S\to C$ be $\P(\mathcal{E})$ over a curve $C$ of positive genus. We can assume $H^0(\mathcal{E})\neq 0$ and $H^0(\mathcal{E}\otimes L)= 0$ for any line bundle $L$ of degree 1. Then $e=-\deg (\mathcal{E})$ is an invariant of $S$. There is a section $C_0$ such that $C_0|_{C_0}=\bigwedge ^2 \mathcal{E}$, which we denote as $D$. Then $\deg(D)=-e$ and 
$$K_S\sim -2C_0+\pi^*(K_C+D).$$ 
Then we have, 
$$0=(K_S+B)\cdot C_0=e+2g-2+tC_0^2+\sum_{B_i\neq C_0}b_i B_i\cdot C_0 \ge 2g-2+ (1-t)e,$$
where $t<1$ is the coefficient of $C_0$ in $B$. So $e \le 0$, and if $e=0$, then $g=1$. On the other hand if $e<0$, (\cite{ha77}, V.2.21) shows that 
$$B\cdot C_0\ge (2C_0+eF)\cdot C_0 \ge -e,$$ so again we conclude $g=1$. And the equality holds only when for each irreducible component $B_i$ of $\supp(B)$, $g_i:B_i\to C$ is unramified.

Now we have $B_i\equiv a(2C_0+eF)$ and $B_i^2=0$.

{\it Claim:} For $g:B_i \to C$, there are line bundles $M$ and $N$ on $B_i$ such that
$$0\to \mathcal{O}\to g^*(\mathcal{E})\otimes M \to N \to 0, \mbox{ and } \deg{N}=0.$$ 

\noindent
{\it Proof of the claim:} it suffices to show that for the induced section $h:B_i\to\P_{B_i}(g_i^*\mathcal{E})$ we have   $h(B_i)^2=0$. But for $B_i\subset \P(\mathcal{E})$,
$$0\to I_{B_i}/I_{B_i}^2\to \Omega^1_{\P(\mathcal{E})}\otimes B_i\to \Omega^1_{B_i}\to 0,$$ 
where $\deg(I_{B_i}/I_{B_i}^2)=0$. Since $B_i\to C$ is unramified, the same exact sequence also computes the conormal bundle of $h(B_i)\subset \P_{B_i}(g_i^*\mathcal{E})$ and hence we conclude that $h(B_i)^2=0$.

If $e>0$, then $\mathcal{E}$ is indecomposable which also means it is stable. 
But the above exact sequence shows that $g^*(E)$ is not semi-stable, which is a contradiction.

When $e=0$, if $\mathcal{E}$ is indecomposable, then $E$ is given by a non-split extension 
$$0\to \mathcal{O} \to \mathcal{E} \to \mathcal{O} \to 0. $$
For any isogeny $g:B\to C$, if there is a surjection $g^*(\mathcal{E})\to F\to 0$ for a line bundle $F$ of degree 0, then it is  given by the pull-back of the above exact sequence. So $B=2C_0$, which contradicts to the assumption that $(S,B)$ is $klt$. Similarly, if $\mathcal{E}=\mathcal{O}\oplus L$ for some degree 0 bundle $L$, then we have at most two candidates for $B_i$, which again contradicts $(S,B)$ being $klt$.

Finally, when $\mathcal{E}=\mathcal{O}^2$, then all candidates of $B_i$ are proportional to $C_0$. Thus the last statement holds.
\end{proof}

\section{Effective Iitaka Fibration}

 In this section we prove the effectiveness of the log Iitaka fibration for  $(X,\Delta)$  a $klt$ pair of dimension three  and four and log  Kodaira codimension two.

Consider a log resolution $\pi:X'\la X$ of $(X,\Delta)$ and write 
$$\pi^*(K_X+\Delta)\equiv K_{X'}+(\pi^{-1})_*\Delta+\sum_ie_iE_i$$ 
with $E_i$ exceptional.
There is a natural number $n$ such that $e_i<1-\frac{1}{n}$ for every $i$. 
Define 
$$\Delta'=(\pi^{-1})_*\Delta+\sum(1-\frac{1}{n})E_i.$$ 
Since all $E_i$ are exceptional divisors, we have  
$$H^0(X',\round{m(K_{X'}+\Delta')})= H^0(X,\round{m(K_X+\Delta)}).$$
By replacing the $\mathcal{A}$ with the DCC set $\mathcal{A}\cup\{1-\frac{1}{n}|n\in\mathbb{N}\}$, we can assume that $X$ is smooth and $\Delta$ is simple normal crossing.
Furthermore, we can assume that there is a morphism $f:X \la C$ giving the Iitaka fibration for $K_X+\Delta$ with $C$ a smooth projective variety of dimension $m$ (cf. \cite{la04}, 2.1.C).  For $F$ the general fiber of $f$, we have that $\kappa(K_F+\Delta_{F})=0$.

From Kawamata's canonical bundle formula, we know that there are a $\Q$-line bundle $M=M^{SS}_{(M,\Delta)/C}$, a $\Q$-divisor $B$ on $C$ and a $\Q$-divisor $R$ on $X$ such that
$$K_X+\Delta\sim_{\Q}f^*(K_C+M+B)+R,$$
where the terms in the formula satisfy the property that if $b(K_F+\Delta_F)$ is Cartier and $h^0(F,b(K_F+\Delta_F))=1$, then for any $r\in \Z_{\ge 0}$ divided by $b$, 
$$H^0(X,\round{r(K_X+\Delta)})=H^0(C,\round{r(K_C+M+B)}).$$
(cf. (\cite{fm}, Theorem 4.5)). $M$ and $B$ are usually referred as the moduli part and the boundary part of the algebraic fiber space. \\

\noindent {\it Proof of (\ref{rel2})}: First we can assume that $(F,\Delta_F)$ is a minimal pair that is  $K_{F}+\Delta_{F}\equiv 0$. In fact, from (\cite{ko07}, Definition 8.4.6), we know the moduli part $M$ only depends on the birational class of the morphism $X\to C$, and we can choose one such that the generic fiber is minimal. Hence Theorem (\ref{epsilon}) implies that $(F,\Delta_{F})$ is $\epsilon=\epsilon(\mathcal{A})$-lc.

 Let $d$ be the smallest positive integer such that $d(K_{F}+\Delta_{F})$ is Cartier and $h^0(F,d(K_{F}+\Delta_{F}))=1$. Then we can construct the corresponding cyclic cover $g:E\to F$, and let $\phi:E'\to E$ be the minimal resolution of $E$.
Then as in (\cite{fm}, Theorem 3.1), it suffices to bound the second Betti number of $E'$.
Write $g^*(K_{F}+\Delta_{F})=K_{E}+\Delta_{E}$. 
From (\cite{km98}, 5.20.3), we know $(E,\Delta_E)$ is $\epsilon$-log canonical, provided that so is $(F,\Delta_F)$, which is the case because of Theorem (\ref{epsilon}). Pulling back to $E'$, we conclude that $(E',\Delta_{E'})$ is $\epsilon$-log canonical. Then (\cite{am04}, Theorem 1.8) implies that $\rho(E')\le 128/\epsilon^5$.  Thus the second Betti number of $E'$ is bounded  by a constant depending only on $\mathcal{A}$.

The remaining part of the proof was given in \cite[Section 3]{tod}. We give a sketch here.
First we assume that the generic fiber $(F,\Delta_F)$ is log smooth, then as in \cite[8.4.5(7)]{ko07}, we can define a local system $\V$ on an open set $U\subset X$. Furthermore, if the the monodromy action on $R^2f_*(\V)$ is unipotent, then the moduli part $M$ is identified with the bottom Hodge filtration of $R^2f_*(\V)$, which is an integral divisor. In general, we take a log resolution $\tilde{F}\to F$ and write 
$$\pi^*(K_F+\Delta_F)=K_{\tilde{F}}+\Delta_{\tilde{F}}-G,$$
where $\round{\Delta_{\tilde{F}}}$, $G\ge 0$ is an integral divisor and $\supp(\Delta_{\tilde{F}})$ has no common component with $G$. Near every codimension one point $P$, there is a Galois cover $Y'$ of $Y$ such that when we pull back everything to $Y'$, the monodromy action becomes unipotent near $P$. Thus as in \cite[3.6]{fm}, the remaining question is to bound the character of $G\to \C^*$ corresponding to the representation of $G$ on $H^0(\tilde{F},\mathcal{O}(G))$. If we write $\tilde{E}$ as the degree-$b$ cyclic cover corresponding to $b(K_{\tilde{F}}+\Delta_{\tilde{F}}-G)\sim 0$, then $H^0(\tilde{F},\mathcal{O}(G))$ is a direct summand of the bottom piece of the Hodge filtration of  $H^2(\tilde{E},\C)$, which is also the bottom piece of $H^2(E',\C)$ since $\tilde{E}$ and $E'$ are birational.

 Then from \cite[3.8]{fm} we conclude that if the index of $(F,\Delta_F)$ is bounded by $b=b(\mathcal{A})$ and the second Betti number of the cyclic cover is bounded by $B=B(\mathcal{A})$, then the denominator of the moduli part is bounded by 
$$a=b\cdot\mbox{lcm}\{m\in\Z_{> 0}|\varphi(m)\le B\},$$
where $\varphi$ denotes the Euler $\varphi$-function.

\begin{thm}\label{iitaka}
Notation as above. Assume that the coefficients of $\Delta$ are in a 
DCC set of rational numbers $\mathcal{A}\subset [0,1]$.  Let the dimension of $X$ be either three or four and the dimension of the general fiber be two. Then there is an 
 constant $N$ depending only on 
the set $\mathcal{A}$  such that
$\round{ N(K_X+\Delta) }$ induces the Iitaka fibration.
\end{thm}

We prove the statement when the dimension of $X$ is four, and leave the easier dimension three case  to the reader. In fact, this follows directly from the argument in (\cite{fm}, Section 6) and Theorem (\ref{rel2}).

\begin{rmk} The main theorem of \cite{vz} says that for a smooth variety $X$ of n arbitrary dimension of  Kodaira dimension 2 and  assuming that the generic fiber of the Iitaka fibration is a smooth variety $F$, there is a constant $N$ depending on the middle Betti number $b_m$ of $F$ and the index of $F$, such that $|NK_X|$ gives the Iitaka fibration.
From (\ref{uniform}) and (\ref{rel2}), we prove that, the analogous bounds for the Betti number and the index, which only depends on the DCC set $\mathcal{A}$, exist even in the log case, provided the generic fiber is of dimension 2. 
\end{rmk} 

\begin{thm}\label{base}
 If $(W,D)$ is a $klt$ surface, $L$ is a nef $\Q$-divisor (not necessarily effective) such that and $K_W+D+L$ is big. Assume that $a$ is a positive integer such that $aL$ is a Cartier divisor, the coefficients of $D$ are in a DCC set $\mathcal{B}$, then we have a uniform $N=N(a,\mathcal{B})$ such that $|\round{N(K_W+D+L)}|$ gives a birational map. 
\end{thm}

{\it Proof of (\ref{iitaka})}: We apply (\ref{base}) to $W=C$, $D=B$, $L=M$. To check the assumptions, the coefficients of $B$ are of the form 
$$\frac{b+n-1}{n}, \mbox{ for some }b\in \mathcal{A}\cup\{0\} \mbox{ and } n\in \Z_{> 0},$$
(cf. \cite{ko07}, Theorem 8.3.7(2)) which forms a DCC set depending on $\mathcal{A}$. And for  $a$, it can be chosen as in the last part of the proof of (\ref{rel2}).\qed

{\it Proof of (\ref{base})}: This is essentially proved in \cite{vz} so here we only give an outline of the proof following \cite{am04}.  As before, we can start by assuming that $W$ is smooth. The main observation is that in fact many  of the  results in \cite{am04} can be strengthened in the way that instead of assuming $(X,B=\sum b_jB_j)$ ($b_j\in \mathcal{B}$) is big, we can assume that $(X,B+L)$ is big, where $L$ is a nef line bundle such $L\cdot C$ is in another DCC set $\mathcal{C}$, then we get the same conclusion by changing all our constants $c=c(\mathcal{B})$ to $c=c(\mathcal{B},\mathcal{C} )$. In particular, we have the following:

\begin{prop}\label{neftail}
Let the notation be as in (\ref{base}). Then there is a uniform $\beta=\beta(\mathcal{B},a)$ such that $K_X+\sum(1-x_j)b_jB_j+L$ is big, provided $x_j\le \beta$. 
\end{prop}

Because of (\ref{neftail}) and the fact that all positive $b_j$ have a lower bound, we can indeed assume that all $b_j$ are of the form $n_j/m$ for some $m=m(\mathcal{B},a)$.

Then the usual argument of cutting log canonical centers (cf. \cite{tod}) works as long as we can prove that the volume of $K_X+\sum(1-x_j)b_jB_j+L$ has a uniform lower bound.

Now we run the minimal model program,  $f:(W,D)\to (W',D'=f_*(D)), L'=f_*L$,  because $f^*(L')\ge L$, we have $$H^0(W,\round{n(K_{W}+D+L)})=H^0(W',\round{n(K_{W'}+D'+L')}), \forall n.$$

\noindent
Case (1): $K_W+D$ is pseudo-effective, then we end with a minimal model $(W',D')$ with $K_{W'}+D'$ is nef. If $K_{W'}+D'$ is big, then the Proposition follows  from \cite{am04}. Otherwise, 
$$(K_{W'}+D'+L)^2=(K_{W'}+D')\cdot L'+L'^2>0,$$
which has a uniform lower bound by our assumption. 

If $K_W+D$ is not pseudo-effective then we define the {\it pseudo-effective threshold}, which is the smallest number $e$ such that $K_W+D+eL$ is pseudo-effective. 

\noindent
Case (2): we end with a log del Pezzo surface $(W',D')$ of Picard number 1. So 
$$K_{W'}+D'+eL'\equiv 0.$$ 

The above discussion about the generalization of \cite{am04} indeed implies that $1-e$ has a uniform lower bound. 
So $(K_{W'}+D'+L')^2=(1-e)L'^2$.

\noindent
Case (3): we end with $(W',D')$ with  a Fano contraction  to a curve. By the above argument we can assume that $L'$ is not big. Then $W'$ is of Picard group of rank 2 generated by the fiber of the Fano contraction and $L'$. Taking the intersection of $K_{W}'+D'+eL'$ with the fiber,  we conclude that $e$ is uniformly far from 1. Since the coefficients of $D'$ have bounded denominators, the fact that $(K_{W'}+D')\cdot L'$ is positive actually implies it is uniformly away from 0. So the volume  
$$(K_{W'}+D'+L')^2=(1-e)(K_{W'}+D')\cdot L'$$
 is bounded from below.

\vspace{5mm}
\noindent Gueorgui Todorov\\
\noindent Department of Mathematics, Princeton University.\\
\noindent {\it Email:} gtodorov@math.princeton.edu\\

\noindent Chenyang Xu\\
\noindent Department of Mathematics, Institute for Advanced Study.\\
\noindent {\it Email:} chenyang@math.ias.edu

\end{document}